\documentclass[11pt,twoside,reqno,centertags,draft]{amsart}
\usepackage{amsfonts}
\usepackage{color,enumitem,graphicx}
\usepackage[colorlinks=true,urlcolor=blue,
citecolor=red,linkcolor=blue,linktocpage,pdfpagelabels,
bookmarksnumbered,bookmarksopen]{hyperref}

  \usepackage{amsmath,amsthm,amsfonts,amssymb}

\begin{document}
\title{ On  semi-linear elliptic equation arising from Micro-Electromechanical Systems  with contacting elastic membrane }
\date{}
\maketitle

\vspace{ -1\baselineskip}

{\small
\begin{center}

\medskip

 {\sc  Huyuan Chen}
  \medskip

 Institute of Mathematical Sciences,  New York University  Shanghai,\\
Shanghai, PR China\\
Email: chenhuyuan@yeah.net\\[10pt]
   {\sc  Ying Wang}
   \medskip

Departamento de Ingenier\'{\i}a  Matem\'atica, Universidad de Chile,\\
  Santiago, Chile\\
 Email: yingwang00@126.com\\[10pt]
 {\sc  Feng Zhou}
   \medskip

Department of Mathematics, East China Normal University,\\
 Shanghai, PR China\\
 Email: fzhou@math.ecnu.edu.cn\\[10pt]

\end{center}
}

\bigskip

\begin{quote}
{\bf Abstract.} This paper is concerned with the nonlinear elliptic problem
 $-\Delta u=\frac{\lambda }{(a-u)^2}$ on a bounded domain $\Omega$ of $\mathbb{R}^N$
 with Dirichlet boundary conditions. This problem arises from Micro-Electromechanical Systems devices
 in the case that the elastic membrane contacts the ground plate on the boundary.
 We analyze  the properties of minimal solutions to this equation when $\lambda>0$
 and the function
  $a:\bar\Omega\to[0,1]$ satisfying
 $a(x)\ge \kappa{\rm dist}(x,\partial\Omega)^\gamma$ for some $\kappa>0$ and $\gamma\in(0,1)$.
Our results show how the boundary decay of the membrane works on the solutions and  pull-in voltage $\lambda$.
\end{quote}

 \renewcommand{\thefootnote}{}

\newcommand{\N}{\mathbb{N}}
\newcommand{\R}{\mathbb{R}}
\newcommand{\Z}{\mathbb{Z}}

\newcommand{\cA}{{\mathcal A}}
\newcommand{\cB}{{\mathcal B}}
\newcommand{\cC}{{\mathcal C}}
\newcommand{\cD}{{\mathcal D}}
\newcommand{\cE}{{\mathcal E}}
\newcommand{\cF}{{\mathcal F}}
\newcommand{\cG}{{\mathcal G}}
\newcommand{\cH}{{\mathcal H}}
\newcommand{\cI}{{\mathcal I}}
\newcommand{\cJ}{{\mathcal J}}
\newcommand{\cK}{{\mathcal K}}
\newcommand{\cL}{{\mathcal L}}
\newcommand{\cM}{{\mathcal M}}
\newcommand{\cN}{{\mathcal N}}
\newcommand{\cO}{{\mathcal O}}
\newcommand{\cP}{{\mathcal P}}
\newcommand{\cQ}{{\mathcal Q}}
\newcommand{\cR}{{\mathcal R}}
\newcommand{\cS}{{\mathcal S}}
\newcommand{\cT}{{\mathcal T}}
\newcommand{\cU}{{\mathcal U}}
\newcommand{\cV}{{\mathcal V}}
\newcommand{\cW}{{\mathcal W}}
\newcommand{\cX}{{\mathcal X}}
\newcommand{\cY}{{\mathcal Y}}
\newcommand{\cZ}{{\mathcal Z}}

\newcommand{\abs}[1]{\lvert#1\rvert}
\newcommand{\xabs}[1]{\left\lvert#1\right\rvert}
\newcommand{\norm}[1]{\lVert#1\rVert}

\newcommand{\loc}{\mathrm{loc}}
\newcommand{\p}{\partial}
\newcommand{\h}{\hskip 5mm}
\newcommand{\ti}{\widetilde}
\newcommand{\D}{\Delta}
\newcommand{\e}{\epsilon}
\newcommand{\bs}{\backslash}
\newcommand{\ep}{\emptyset}
\newcommand{\su}{\subset}
\newcommand{\ds}{\displaystyle}
\newcommand{\ld}{\lambda}
\newcommand{\vp}{\varphi}
\newcommand{\wpp}{W_0^{1,\ p}(\Omega)}
\newcommand{\ino}{\int_\Omega}
\newcommand{\bo}{\overline{\Omega}}
\newcommand{\ccc}{\cC_0^1(\bo)}
\newcommand{\iii}{\opint_{D_1}D_i}

\numberwithin{equation}{section}

\vskip 0.2cm \arraycolsep1.5pt
\newtheorem{lemma}{Lemma}[section]
\newtheorem{theorem}{Theorem}[section]
\newtheorem{definition}{Definition}[section]
\newtheorem{proposition}{Proposition}[section]
\newtheorem{remark}{Remark}[section]
\newtheorem{corollary}{Corollary}[section]

\setcounter{equation}{0}
\section{Introduction}

Micro-Electromechanical Systems (MEMS) are often used to combine electronics with micro-size mechanical devices in the design
of various types of microscopic machinery. They
are successfully utilized in components of many commercial systems, including accelerometers for airbag deployment in automobiles, ink jet printer heads, optical switches, chemical sensors, etc.
In  MEMS devices, a key component is called the electrostatic actuation, which
is based on an electrostatic-controlled tunable, it is a simple idealized electrostatic device. The upper part of this electrostatic device consists of a thin and deformable elastic membrane that is fixed along its
boundary and which lies above a rigid grounded plate. This elastic membrane is modeled as a dielectric with
a small but finite thickness. The upper surface of the membrane is coated with a negligibly thin metallic
conducting film. When a voltage $\lambda$ is applied to the conducting film, the thin dielectric membrane deflects
towards the bottom plate, and when $\lambda$ is increased beyond a certain critical value $\lambda^*-$known as pull-in
voltage$-$the steady state of the elastic membrane is lost, and proceeds to touchdown or snap through at a
finite time creating the so-called pull-in instability.

A mathematical model of the physical phenomena, leading to a partial
differential equation for the dimensionless deflection of the membrane, was
derived and analyzed in \cite{GG,G,GHW,GPW,LY,P,YZ} and reference therein. In the damping-dominated limit, and using a
narrow-gap asymptotic analysis, the dimensionless deflection $u$ of the membrane
on a bounded domain $\Omega$ in $\R^2$ is found to satisfy the
equation
\begin{equation}\label{eq 1.1240}
 -\Delta    u = \frac{\lambda  }{(1-u)^2} \qquad {\rm in}\quad \Omega
\end{equation}
with the Dirichlet boundary condition. Here the term on the right hand side of equation (\ref{eq 1.1240}) is the Coulomb force.
Later on,  Ghoussoub and  Guo in \cite{GG,G} studied the nonlinear elliptic problem
\begin{equation}\label{eq 1.1241}
 -\Delta    u = \frac{\lambda f(x) }{(1-u)^2} \qquad {\rm in}\quad \Omega
\end{equation}
 with the Dirichlet boundary condition, which  models a simple electrostatic
MEMS device consisting of a thin dielectric elastic membrane with boundary supported at $0$ above a
rigid ground plate located at $1$.
 Here $\Omega$ is a bounded domain of $\R^N$ and
 the function $f\ge0$ represents the permittivity profile and $\lambda >0$ is a constant which is increasing with respect to the applied voltage.
 We know that for any given suitable $f$, there exists a critical value
$\lambda^*$  (pull-in voltage) such that if $\lambda\in(0,\lambda^*)$,  problem (\ref{eq 1.1241})
is solvable, while for $\lambda>\lambda^*$, no solution for (\ref{eq 1.1241}) exists.

In an effort to achieve better MEMS design,   the
membrane can be technologically fabricated into  non-flat shape like the surface of a semi-ball, which contacts the ground plate along the boundary.
In this paper, we study how the shape of  membranes effects on the pull-in voltage.
In what follows,
we denote that $\Omega$ is a $C^2$ bounded domain in $\R^N$ with $N\ge 1$, $\rho(x)={\rm dist}(x,\partial\Omega)$ for $x\in\Omega$,
the function  $a:\bar\Omega\to[0,1]$ is in the class of $C^\gamma(\Omega)\cap C(\bar\Omega)$ and satisfies
\begin{equation}\label{1.1}
a(x)\ge \kappa\rho(x)^\gamma,\quad \forall \ x\in\Omega
\end{equation}
for some $\kappa>0$ and $\gamma\in(0,1)$.   Our purpose of this paper is to consider the  minimal solutions to
  elliptic equation
\begin{equation}\label{eq 1.1}
\left\{\arraycolsep=1pt
\begin{array}{lll}
 -\Delta    u = \frac{\lambda }{(a-u)^2}\quad  &{\rm in}\quad\ \Omega,
 \\[2mm]
 \phantom{- }
 0<u<a\quad &{\rm in}\quad\ \Omega,
 \\[2mm]
 \phantom{-\Delta   }
 u=0\quad &{\rm on}\quad   \partial \Omega,
\end{array}
\right.
\end{equation}
where parameter $\lambda>0$ characterizes the relative strength of the electrostatic and mechanical forces in the system. Equation (\ref{eq 1.1}) models a MEMS device
that the static deformation  of the surface of membrane when it is applied  voltage $\lambda$,
where $a$ is   initially undeflected state  of the elastic membrane that contacts the ground plate on the boundary.
For this equation, we have the following existence results.

\begin{theorem}\label{teo 1}
Assume that $a\in C^\gamma(\Omega)\cap C(\bar\Omega)$ satisfies (\ref{1.1}) with $\gamma\in(0,\frac23]$ and $\kappa>0$, then there exists a
finite pull-in voltage $\lambda^*:=\lambda^*(\kappa,\gamma)>0$ such that

$(i)$ for $\lambda\in(0,\lambda^*)$, problem (\ref{eq 1.1}) admits a minimal solution $u_\lambda$ and the mapping: $\lambda\mapsto u_\lambda$
is   strictly  increasing;

$(ii)$  for $\lambda>\lambda^*$, there is no solution for (\ref{eq 1.1});

$(iii)$ assume more that there exists $c_0\ge \kappa$ such that
\begin{equation}\label{a 1.0}
a(x)\le c_0\rho(x)^{\gamma}, \quad x\in\Omega,
\end{equation}
 then there exists  $\lambda_*:=\lambda_*(\kappa,\gamma)\in(0, \lambda^*]$ such that for $\lambda\in(0,\lambda_*)$,   $u_\lambda\in H^1_0(\Omega)$ and
$$
{\rm for}\quad \gamma\not=\frac12,\qquad \frac{\lambda}{c_1}\rho(x)^{\min\{1,2-2\gamma\}}\le u_\lambda(x)\le c_1 \lambda \rho(x)^{\min\{1,2-2\gamma\}},\qquad \forall x\in
\Omega;
$$
$$
{\rm for}\quad\gamma=\frac12,\qquad \frac{\lambda}{c_1}\rho(x)\ln \frac1{\rho(x)} \le u_\lambda(x)\le c_1 \lambda \rho(x)\ln \frac1{\rho(x)},\qquad \forall x\in A_{\frac12},
$$
where $c_1\ge 1$ and $A_{\frac12}=\{x\in\Omega:\ \rho(x)<\frac12\}$.

In particular, assume more that  $\Omega=B_1(0)$ and
$$a(x)=\kappa(1-|x|^2)^\gamma,\qquad \forall x\in B_1(0),$$
then the mappings: $\gamma\mapsto \lambda_*(\kappa,\gamma)$, $\gamma\mapsto\lambda^*(\kappa,\gamma)$ are decreasing and
$$
(2/3)^5\kappa^3\le \lambda_*(\kappa,2/3)\le \lambda_*(\kappa,\gamma)\le \lambda^*(\kappa,\gamma)\le \lambda^*(\kappa,0)\le (4/3)^2N\kappa^3.
$$
\end{theorem}

We remark that  the membrane contacts the ground plate on the boundary with decay rate $\rho^\gamma$, $\gamma\in(0,\frac23]$,
there still has a positive finite pull-in voltage $\lambda^*$, but the decay of $a$ plays an important role in decay of  minimal solution,
the regularity of minimal solution and the estimate of $\lambda^*$.   Theorem \ref{teo 1} shows that  the membrane  of the MEMS device could be
designed as the surface of the unit semi-ball,
that is,
$$\Omega=B_1(0)\quad{\rm and}\quad a(x)=(1-|x|^2)^\frac12,$$
which is equivalent to the case that $a(x)=\rho(x)^{\frac12}$, so there exists a positive finite pull-in voltage $\lambda^*$.
However,
the decay rate of function $a$ determined completely
nonexistence of pull-in voltage when $\gamma>\frac23$. Precisely, we have the following non-existence result.

\begin{theorem}\label{teo 2}
Assume that $a\in  C(\bar\Omega)$ is positive and satisfies (\ref{a 1.0}) with $\gamma\in(\frac23,1)$ and $c_0>0$.
Then problem
(\ref{eq 1.1}) admits no nonnegative solution for any $\lambda>0.$

\end{theorem}

We notice that for $\gamma\le\frac23$  and fixed $\kappa$, the finite pull-in voltage $\lambda^*$ depends on $\gamma$, however,
when $\gamma=\frac23$, $\lambda^*>0$ and $\lambda^*=0$ for $\gamma>\frac23$. Therefore, there is a gap at $\gamma=\frac23$. From Theorem \ref{teo 2},
we learn that the   membrane  of the MEMS device should be
made steeply enough near the boundary, otherwise the device does not work.

\smallskip

As normal, it is  challenging to study the extremal solution, i.e., the solution of (\ref{eq 1.1}) when $\lambda=\lambda^*$. Especially, the decay of function $a$
makes this issue more subtle. From Theorem \ref{teo 1}, we observe that the mapping $\lambda\mapsto u_\lambda$ is increasing and
uniformly bounded by function $a$, then it is   well-defined that
\begin{equation}\label{4.2}
 u_{\lambda^*}:=\lim_{\lambda\to\lambda^*} u_\lambda\quad{\rm in}\ \ \bar\Omega,
\end{equation}
where $u_\lambda$ is the minimal solution of (\ref{eq 1.1}) with $\lambda\in(0,\lambda^*)$.
  Our final purpose in this paper is to
  prove that $u_{\lambda^*}$ is a solution of (\ref{eq 1.1}) in some weak sense and it is called the extremal solution. The extremal solution always is found in a weak sense and then it could be improved the regularity until
  to the classical sense when the dimension $N$ is suitable. 

\begin{definition}\label{def 1}
A function $u$ is a weak solution of (\ref{eq 1.1}) if $0\le u\le a$ and
$$
\int_\Omega u (-\Delta \xi) dx=\int_\Omega \frac{\lambda \xi}{(a-u)^2} dx,\quad \forall \xi\in C_c^2(\Omega),
$$
where  $C_c^2(\Omega)$ is the space of all  $C^2$ functions with compact support in $\Omega$.

A solution (or weak solution) $u$ of (\ref{eq 1.1}) is stable (resp. semi-stable) if
$$
\int_\Omega  |\nabla\varphi|^2 dx>\int_\Omega \frac{2\lambda \varphi^2}{(a-u)^3} dx,\quad ({\rm resp.}\ \ge)\quad \forall \varphi\in H^1_0(\Omega)\setminus\{0\}.
$$

\end{definition}

We make use of the functions' space $C^2_c(\Omega)$ in  Definition \ref{def 1}  replacing $C_0^2(\Omega):=C^2(\Omega)\cap C_0(\Omega)$, used in the case of $a\equiv1$, to avoid the singularity at the boundary of $\frac{1}{(a-u_{\lambda})^2}$.

\begin{theorem}\label{teo 3}
Assume that $\lambda\in(0,\lambda^*)$, the function
$a\in C^\gamma(\Omega)\cap C(\bar\Omega)$  satisfies (\ref{1.1}) and(\ref{a 1.0}) with
 $c_0\ge \kappa>0$, $\gamma\in(0,\frac23]$,
$u_\lambda$ is the minimal solution of (\ref{eq 1.1}) and $ u_{\lambda^*}$ is given by (\ref{4.2}).
Then

 $(i)$ \ \ $ u_{\lambda^*}$ is a weak solution of (\ref{eq 1.1}) and $u_{\lambda^*}\in W^{1,\frac{N}{N-\beta}}_0(\Omega)$
for any $\beta\in(0,\gamma)$.

$(ii)$ $u_\lambda$ is a stable solution of (\ref{eq 1.1}) with $\lambda\in(0,\lambda_*)$.

$(iii)$ if $\gamma=\frac23$, we have that $\lambda^*=\lambda_*$,
  $ u_{\lambda^*}$ is a semi-stable weak  solution of (\ref{eq 1.1}).\smallskip

Assume more that $1\le N\le 7$,  $\Omega=B_1(0)$ and $a(x)=\kappa(1-|x|^2)^{\frac23}$, $ u_{\lambda^*}$ is a classical solution of (\ref{eq 1.1}).

\end{theorem}

We note that $(i)$ the extremal solution $u_{\lambda^*}$ is a weak solution but in a
weaker sense comparing with the case $a\equiv1$;
$(ii)$ the main difficulty  to study
 the stability of $u_\lambda$ with $\lambda\in(0,\lambda_*)$ arises from the decay on the boundary of $a$ and to overcome this difficult,
we make use of the generalized Hardy's inequality;
$(iii)$  in the case of $\gamma=\frac23$, the regularity of the extremal solution could be improved into classical sense
 for $1\le N\le 7$,  the same range of the dimension with the case $a\equiv1$.

The paper is organized as follows.  Section 2 is devoted to do the estimate for $\mathbb{G}_\Omega[\rho^{\tau-2}]$ when $\tau\in(0,2)$ and
then we prove that there exists pull-in voltage $\lambda^*$ such that problem (\ref{eq 1.1}) admits a minimal solution when $\lambda\in(0,\lambda^*)$. We also analyze
 the boundary decay of the minimal solution.
 In Section  3, we do the estimate for $\lambda_*$ and $\lambda^*$ when $\Omega=B_1(0)$ and $a(x)=\kappa(1-|x|)^\gamma$. In
Section 4, we study the properties of $u_{\lambda^*}$, including its regularity and   stability.

\setcounter{equation}{0}
\section{Existence }
Denote by $G_\Omega$ the Green kernel of $-\Delta$ in $\Omega\times\Omega$ and by $\mathbb{G}_\Omega[\cdot]$ the
Green operator defined as
$$\mathbb{G}_\Omega[f](x)=\int_{\Omega} G_\Omega(x,y)f(y)dy,\quad \forall f\in L^1(\Omega,\rho^{-1}dx). $$
To do the existence of  a minimal solution of problem (\ref{eq 1.1}), the following estimates play an important role.
\begin{lemma}\label{lm 2.1}
Let $\tau\in(0,2)$, $A_{\frac12}=\{x\in\Omega: \ \rho(x)<\frac12\}$.  For $x\in A_{\frac12}$, we denote
\begin{equation}\label{varrho}
\varrho_\tau(x)=\left\{ \arraycolsep=1pt
\begin{array}{lll}
 \rho(x)^{\min\{1,\tau\}}\ \ &{\rm if}\ \ \tau\in (0,1)\cup(1,2),
 \\[2mm]
\rho(x) \ln \frac1{\rho(x)} \ & {\rm if}\ \ \tau=1
\end{array}
\right.
\end{equation}
and we make $C^1$ extension of $\varrho_\tau$ into  $\Omega\setminus A_{\frac12}$ such that
$\varrho_\tau>0$ in $\Omega\setminus A_{\frac12}$.

Then there exists $c_\tau>1$ such that
$$\frac1c_\tau\varrho_\tau(x) \le \mathbb{G}_\Omega[\rho^{\tau-2}](x)\le c_\tau\varrho_\tau(x),\quad\forall x\in  \Omega.$$

\end{lemma}
{\it Proof.} Since the domain $\Omega$ is  $C^2$, then there exists $\delta_1>0$  such that the distance function
$\rho(\cdot)$ is  $C^2$ in
$A_{\delta_1}=\{x\in \Omega:\ \rho(x)<\delta_1\}.$
Let us define
\begin{equation}\label{2.6}
V_\tau(x)=\left\{ \arraycolsep=1pt
\begin{array}{lll}
 l(x),\ \ \ \ \ &
x\in \Omega\setminus A_{\delta_1},\\[2mm]
 \rho(x)^{\tau},\ \  & x\in A_{\delta_1},
\end{array}
\right.
\end{equation}
 where  $\tau$ is a parameter in $(0,2)$ and  $l$  is a positive function such that
$V_\tau$ is $C^2$ in $\Omega$.
Now we do the estimate of $-\Delta V_\tau$ near the boundary of $\Omega$.   Without loss of generality, we may assume that $0\in\partial\Omega$ and
$e_N=(0,\cdots,0,1)$ is the unit normal vector at $0$ pointing inside. We observe that there exists $\delta_2\in(0,\delta_1)$ such that
for any $t\in(0,\delta_2)$,  $te_N\in A_{\delta_1}$
and  $V_\tau(te_N)=\rho(te_N)^{\tau}=t^\tau$. Then
$$\frac{\partial^2 V_\tau(te_N)}{\partial x_N^2}=\tau(\tau-1)t^{\tau-2} \quad {\rm and}\quad|\frac{\partial^2 V_\tau(te_N)}{\partial x_i^2}|\le c_2,\quad i=1,2,\cdots,N-1,$$
where $c_2>0$ is independent of $t$.
Thus,
$$\tau(\tau-1)t^{\tau-2}-c_2\le\Delta V_\tau(t e_N)\le \tau(\tau-1)t^{\tau-2}+c_2.$$
Since $\Omega$ is a $C^2$ bounded domain, $\partial\Omega$ is compact and then
\begin{equation}\label{2.1.1}
\tau(\tau-1)\rho(x)^{\tau-2}-c_2\le \Delta V_\tau(x)\le \tau(\tau-1)\rho(x)^{\tau-2}+c_2,\quad \forall x\in A_{\delta_2}.
\end{equation}

For $\tau\in(0,1)$, we observe that  $\tau(\tau-1)<0$ and by (\ref{2.1.1}),
$$
\frac{1}{2}\tau(1-\tau)\rho(x)^{\tau-2}\le -\Delta V_\tau(x)\le 2\tau(1-\tau)\rho(x)^{\tau-2},\quad \forall x\in A_{\delta_2}.
$$
Combining with the fact that
$-\Delta\mathbb{G}_\Omega[\rho^{\tau-2}]=\rho^{\tau-2}$ in $\Omega$,
we obtain that
$$\frac{1}{2\tau(1-\tau)}(-\Delta V_\tau)\le-\Delta\mathbb{G}_\Omega[\rho^{\tau-2}]\le \frac{2}{\tau(1-\tau)}(-\Delta V_\tau)
\quad {\rm in}\ \ A_{\delta_2}.
$$
For $x\in\partial\Omega$, we have that $V_\tau(x)=0$ and $G_\Omega(x,y)=0$ for $y\in\Omega$ and then
$$\mathbb{G}_\Omega[\rho^{\tau-2}](x)=\int_{\Omega} G_\Omega(x,y)\rho(y)^{\tau-2}dy=0.$$
  For $x\in \Omega$ satisfying $\rho(x)=\delta_2$, we have that $V_\tau(x)={\delta_2}^\tau$ and then
\begin{equation}\label{2.1.3251}
\frac1{c_3}V_{\tau}(x)\le \mathbb{G}_\Omega[\rho^{\tau-2}](x)\le c_3V_{ \tau}(x)\quad{\rm for}\ \rho(x)=\delta_2,
\end{equation}
where $c_3>1$.
By Comparison Principle,   there exists $c_4>1$ depends on $\tau$ such that
\begin{equation}\label{2.1.3252}
\frac1{c_4}V_{\tau}\le \mathbb{G}_\Omega[\rho^{\tau-2}]\le c_4V_{ \tau}\quad{\rm in}\ \  A_{\delta_2}.
\end{equation}
Since $\mathbb{G}_\Omega[\rho^{\tau-2}]$ and $V_{\tau}$ is bounded in $\Omega\setminus A_{\delta_2}$, then (\ref{2.1.3252})
holds in $\Omega$.

For $\tau\in(1,2)$,
by the fact of (\ref{2.1.1}), we observe that
$$
\frac{1}{2}\tau(\tau-1)\rho(x)^{\tau-2}\le \Delta V_\tau(x)\le 2\tau(\tau-1)\rho(x)^{\tau-2},\quad \forall x\in A_{\delta_2}.
$$
Letting $W_\tau=\mathbb{G}_\Omega[1]-V_\tau$, then there exists  $c_5>1$ depends on $\tau$ such that
$$\frac{1}{c_5}\rho(x)^{\tau-2}\le -\Delta W_\tau(x)=-1+ \Delta V_\tau(x)\le  c_5\rho(x)^{\tau-2},\quad \forall x\in A_{\delta_2}.$$
By the fact that
$-\Delta\mathbb{G}_\Omega[\rho^{\tau-2}]=\rho^{\tau-2}$ in $\Omega$, we have
\begin{equation}\label{2.1.3253}
\frac1{c_6}W_{\tau}(x)\le \mathbb{G}_\Omega[\rho^{\tau-2}](x)\le c_6W_{ \tau}(x), \ \ \ \ \forall x\in\partial A_{\delta_2}
\end{equation}
for some $c_6>1$ depends on $\tau$.
Using the Comparison Principle, there exists $c_7,c_8>1$ depend on $\tau$ such that
$$\frac1{c_8}\rho(x)\le \frac1{c_7}W_\tau(x) \le \mathbb{G}_\Omega[\rho^{\tau-2}](x)\le c_7W_\tau(x) \le c_8\rho(x), \quad \  x\in A_{\delta_2}, $$
and then
$$\frac1{c_8}\rho(x)\le \mathbb{G}_\Omega[\rho^{\tau-2}](x)\le c_8\rho(x), \quad \  x\in \Omega. $$

 For $\tau=1$, we redefine
$$
 V_1(x)=\left\{ \arraycolsep=1pt
\begin{array}{lll}
 l(x),\ \ \ \ \ &
x\in \Omega\setminus A_{\delta_1},\\[2mm]
\rho(x)\ln\frac1{\rho(x)},\ \ \ & x\in A_{\delta_1},
\end{array}
\right.
$$
 where   $l$  is a positive function such that
$V_1$ is $C^2$ in $\Omega$. By direct computation, there exist $\delta_3\in(0,\delta_1)$ and $c_9>0$ such that
$$
 \rho(x)^{ -1}-c_9\le -\Delta V_1(x)\le  \rho(x)^{ -1}+c_9,\quad \forall x\in A_{\delta_3}.
$$
Then it follows  by Comparison Principle that
$$\frac1{c_{10}} \rho\ln \frac1\rho \le \mathbb{G}_\Omega[\rho^{-1}]\le  c_{10} \rho\ln \frac1\rho \qquad{\rm in }\quad  A_{\delta_3} $$
for some $c_{10}>1$.
The proof is complete. \qquad $\Box$

\smallskip

By Lemma \ref{lm 2.1}, we have the following results.
\begin{corollary}\label{cr 2.1}
 For $\gamma\in(\frac23,1)$, we have that
\begin{equation}\label{2.1.3}
 \lim_{x\in\Omega, x\to\partial\Omega}\mathbb{G}_\Omega[\rho^{-2\gamma}](x)\rho(x)^{-\gamma}=+\infty.
\end{equation}
\end{corollary}
{\it Proof.} Let $\tau=2-2\gamma$, by $\gamma\in(\frac23,1)$, we have that
$\tau<\gamma<1$. Using Lemma \ref{lm 2.1}, we obtain that
$$\mathbb{G}_\Omega[\rho^{-2\gamma}](x)\ge c_{11}\rho(x)^{\tau},\quad \forall x\in\Omega$$
for some $c_{11}>0$, which implies (\ref{2.1.3}).
\qquad $\Box$

\smallskip

Now we are ready to show the existence of pull-in voltage $\lambda^*$ to problem (\ref{eq 1.1}) such that
 (\ref{eq 1.1}) admits a solution for $\lambda\in(0,\lambda^*)$ and no solution for $\lambda>\lambda^*$.

\begin{proposition}\label{pr 2.1}
 Assume that $a\in C^\gamma(\Omega)\cap C(\bar\Omega)$ satisfies (\ref{1.1}) with $\gamma\in(0,\frac23]$, then there exists  $\lambda^* >0$ such that
 problem  (\ref{eq 1.1}) admits at least one solution for  $\lambda\in(0,\lambda^*)$ and no solution
for $\lambda>\lambda^*$.
Moreover,
\begin{equation}\label{2.4}
\lambda^*\le \frac{\int_\Omega a(x) dx}{ \int_\Omega\frac{\mathbb{G}_\Omega[1](x)}{a(x)^2}dx}.
\end{equation}

\end{proposition}
{\it Proof.}   Let $v_0\equiv0$ in $\bar\Omega$ and
$$v_1=\lambda \mathbb{G}_\Omega[\frac1{a^2}]>0,$$
by (\ref{1.1}) and Lemma \ref{lm 2.1},
\begin{equation}\label{2.428}
 v_1  = \lambda \mathbb{G}_\Omega[\frac1{a^2}] \le  \frac{\lambda}{\kappa^2} \mathbb{G}_\Omega[\rho^{-2\gamma}]  \le  \frac{ \lambda}{\kappa^2}c_{12} \varrho_{ 2-2\gamma },
\end{equation}
where $c_{12}>0$ depending on $\gamma$ and $\varrho_{ 2-2\gamma }$ is given by (\ref{varrho}).

For $0<\gamma\le \frac23$ and $\gamma\not=\frac12$, we observe that $\min\{1,2-2\gamma\}\ge \gamma$ and by (\ref{2.428}),
$$
v_1(x)\le  \frac{ \lambda}{\kappa^2}c_{12} \rho(x)^{\min\{1,2-2\gamma\}}\le \frac{\lambda}{\kappa^2}c_{13}\rho(x)^\gamma,\qquad x\in\Omega.
$$
For $\gamma=\frac12$, we see that $2-2\gamma=1$ and by (\ref{2.428}),
$$
v_1(x)\le  \frac{ \lambda}{\kappa^2}c_{12} \rho(x)\log\frac1{\rho(x)}\le \frac{\lambda}{\kappa^2}c_{13}\rho(x)^\gamma,\qquad x\in\Omega.
$$
Then
$$v_1(x)\le \frac{\lambda}{\kappa^2}c_{13}\rho(x)^\gamma,\qquad x\in\Omega.$$
Fixed $\mu\in(0, \kappa)$, there exists $\lambda_1>0$ such that
$$\frac{\lambda}{\kappa^2}c_{13}\le \mu<\kappa, \ \ \ \forall\lambda<\lambda_1,$$
then  for any $\lambda<\lambda_1$,
$$v_1(x)\le \mu\rho(x)^\gamma,\qquad \forall x\in\Omega.$$

Let $v_2=\lambda \mathbb{G}_\Omega[\frac1{(a-v_1)^2}]$, by the fact that $a(x)\ge \kappa\rho(x)^\gamma>\mu\rho(x)^\gamma\ge v_1(x)>0$ for $x\in\Omega$ and and Lemma \ref{lm 2.1},
 we have that
$$v_1=\lambda\mathbb{G}_\Omega[\frac1{a^2}]\le v_2\le  \frac{\lambda}{(\kappa-\mu)^2}\mathbb{G}_\Omega[\rho^{-2\gamma}]
\le   \frac{\lambda}{(\kappa-\mu)^2}c_{12} \varrho_{ 2-2\gamma }\le  \frac{\lambda}{(\kappa-\mu)^2}c_{14}\rho^\gamma\quad {\rm in}\ \ \Omega$$
for $\gamma\in(0,\frac23]$. We observe that there exists $\lambda_2\in(0,\lambda_1]$ such that
$$\frac{\lambda}{(\kappa-\mu)^2}c_{14}\le \mu<\kappa, \ \ \ \forall\lambda<\lambda_2.$$
Then for any $\lambda<\lambda_2$, we obtain that
$$v_2(x)\le \mu\rho(x)^\gamma,\qquad \forall x\in\Omega.$$

Iterating above process, we have that
 $$v_n:=\lambda\mathbb{G}_\Omega[\frac{1}{(a-v_{n-1})^2}]\ge v_{n-1},\quad n\in\N$$
and
$$
v_n(x)\le \mu\rho(x)^\gamma,\qquad\forall x\in\Omega.
$$
Thus, the sequence $\{v_n\}_n$ converges, denoting by $u_\lambda=\lim_{n\to\infty} v_n$, then $u_\lambda$ is  a classical solution of (\ref{eq 1.1}).

We claim that $u_\lambda$ is the minimal solution of  (\ref{eq 1.1}), that is, for any positive solution $u$ of (\ref{eq 1.1}), we always have
$u_\lambda\le u$. Indeed, there holds $u\ge v_0$ and then
$$u=\lambda \mathbb{G}_\Omega[\frac1{(a-u)^2}]\ge\lambda\mathbb{G}_\Omega[\frac1{a^2}]=v_1.$$
We may show inductively that $$u\ge v_n$$
for all $n\in\N$. The claim follows.

Similarly, if problem  (\ref{eq 1.1})
 has a  super solution $u$  for $\lambda_0>0$,
then (\ref{eq 1.1}) admits a minimal solution $u_\lambda$   for all $\lambda\in(0,\lambda_0]$. As a consequence,
   the mapping $\lambda\mapsto u_\lambda$ is increasing. So we may define
   $$\lambda^*=\sup\{\lambda>0:\,  (\ref{eq 1.1}) \ {\rm has \  a \  minimal \  solution \ for} \ \lambda \},$$
which is the largest $\lambda$ such that problem (\ref{eq 1.1}) has minimal positive solution, and $\lambda^*>0$.

 For $0<\lambda_1<\lambda_2<\lambda^*$, we know that $0\le u_{\lambda_1}\le u_{\lambda_2}\le a$ in $\Omega$, then
\begin{eqnarray*}
 -\Delta (u_{\lambda_2}-u_{\lambda_1}) &=& \frac{\lambda_2}{(a-u_{\lambda_2} )^2} - \frac{\lambda_1}{(a-u_{\lambda_1} )^2}\\
    &\ge & \frac{\lambda_2-\lambda_1}{(a-u_{\lambda_1} )^2}\ge \frac{\lambda_2-\lambda_1}{a^2}>0,
\end{eqnarray*}
which implies that
 \begin{equation}\label{e 1.1}
 u_{\lambda_2}-u_{\lambda_1}\ge (\lambda_2-\lambda_1)\mathbb{G}_\Omega[a^{-2}]>0.
 \end{equation}
It infers that $u_{\lambda}<a$ in $\Omega$ for any $\lambda<\lambda^*$ and
by the interior regularity, we have that $u_\lambda\in C_{loc}^{2,\gamma'}(\Omega)$ for any $\gamma'<\gamma$.

Now we prove that $\lambda^*<+\infty$.  If not, then for any $\lambda>0$, problem  (\ref{eq 1.1}) has the minimal solution $u_\lambda$.
Let $A_{\delta}=\{x\in\Omega:\, \rho(x)<\delta\}$ and $n\in\N$,
\begin{equation}\label{101}
\eta_n=1 \ \ {\rm in} \  \Omega\setminus{A_{1/n}}, \quad  \   \eta_n=0 \ \ {\rm in} \  A_{1/{2n}},\quad  \  \eta_n\in C^2(\Omega),
\end{equation}
and
$\xi_n= \mathbb{G}_{\Omega}[1]\eta_n$ in $\Omega$, we observe that $\xi_n\in C_c^2(\Omega)$ and
\begin{eqnarray*}
\int_\Omega u_\lambda (-\Delta \xi_n) dx=\int_\Omega u_\lambda (\eta_n-\nabla \mathbb{G}_{\Omega}[1]\cdot\nabla\eta_n+(-\Delta \eta_n)\mathbb{G}_{\Omega}[1])dx.
\end{eqnarray*}
By the fact that $|\nabla \eta_n|\le c_{15}n$, $|(-\Delta)\eta_n|\le c_{16}n^2$ and
$\mathbb{G}_{\Omega}[1]\le c_{17}\rho$ in $\Omega$, then we have that
\begin{equation}\label{eq9281}
\int_\Omega u_\lambda (-\Delta\xi_n) dx\le\int_\Omega u_\lambda dx+c_{15}n\int_{A_{1/n}}u_\lambda dx + c_{18}n^2\int_{A_{1/n}}u_\lambda \rho dx
\le c_{19},
\end{equation}
where $c_{15}, c_{16}, c_{17}, c_{18}, c_{19}>0$ are independent on $n$.

Since $u_\lambda$ is the minimal solution of  (\ref{eq 1.1}), then
$$\int_\Omega u_\lambda (-\Delta\xi_n) dx=\int_\Omega \nabla u_\lambda\cdot \nabla{\xi_n} dx
=\int_\Omega (-\Delta u_\lambda)\xi_n dx=\int_\Omega \frac{\lambda \xi_n}{(a-u_\lambda)^2} dx$$
Passing to the limit of $n\to\infty$ and combining with (\ref{eq9281}), we see that
$$\int_\Omega \frac{\lambda \mathbb{G}_{\Omega}[1]}{(a-u_\lambda)^2} dx\le c_{20}.$$
Thus,
\begin{eqnarray*}
\int_\Omega a(x) dx\ge \int_\Omega u_\lambda(x) dx&=&\int_\Omega \mathbb{G}_\Omega[1](x) (-\Delta u_\lambda) (x) dx\\
 &=& \lambda \int_\Omega \frac{\mathbb{G}_\Omega[1](x)}{[a(x)-u_\lambda(x)]^2}dx
  \ge  \lambda \int_\Omega\frac{\mathbb{G}_\Omega[1](x)}{a^2(x)}dx,
\end{eqnarray*}
which implies that
$$
\lambda\le  \frac{\int_\Omega a(x) dx }{\int_\Omega\frac{\mathbb{G}_\Omega[1](x)}{a^2(x)}dx}.
$$
Therefore,
$$
\lambda^*\le  \frac{\int_\Omega a(x) dx }{\int_\Omega\frac{\mathbb{G}_\Omega[1](x)}{a^2(x)}dx}.
$$

Similar proof as the claim, we know that  problem  (\ref{eq 1.1}) has the minimal solution for  $\lambda\in(0,\lambda^*)$ and no solution
for $\lambda>\lambda^*$.
This end the proof.
\qquad $\Box$
\medskip

\noindent{\bf Proof of Theorem \ref{teo 2}.}
By contradiction. Under the assumptions of Theorem \ref{teo 2}, if there exists some $\lambda>0$ such that
problem (\ref{eq 1.1}) has a solution $u_\lambda$ satisfying $0<u_\lambda<a$ in $\Omega$, then
\begin{equation}\label{2.5}
 \lambda \mathbb{G}_\Omega[\frac1{a^2}]\le \lambda \mathbb{G}_\Omega[\frac1{(a-u_\lambda)^2}] = u_\lambda< a\quad{\rm in}\quad \Omega.
\end{equation}
On the other hand, by (\ref{a 1.0}), we have that
$$  \lambda \mathbb{G}_\Omega[\frac1{a^2}]\ge \frac{\lambda}{c_0^2}\mathbb{G}_\Omega[\rho^{-2\gamma}].$$
For $\gamma\in(\frac23,1)$, using Corollary \ref{cr 2.1} and (\ref{2.5}), we deduce that
$$\lim_{x\in\Omega,  x\to\partial\Omega} a(x)\rho(x)^{-\gamma}=+\infty,$$
which contradicts (\ref{a 1.0}). This ends the  proof.\qquad$\Box$

\smallskip

We next do the boundary decay estimate for $u_\lambda$. To this end, we introduce  the following lemma.
\begin{lemma} \label{lm 2.2}

Assume that the function $a\in C^\gamma(\Omega)\cap C(\bar\Omega)$  satisfies (\ref{1.1}) with  $\gamma\in(0,\frac23]$ and  $u$ is a super solution of (\ref{eq 1.1})
with $\lambda>0$ such that
\begin{equation}\label{2.3}
u \le \theta a\quad {\rm in}\quad \Omega
\end{equation}
for  some $\theta\in(0,1)$,
then (\ref{eq 1.1}) admits the minimal solution $u_\lambda$  such that
\begin{equation}\label{2.311} u_\lambda\le c_{20}\varrho_{ 2-2\gamma }\quad{\rm in}\quad \Omega\end{equation}
 for some $c_{20}>0$, where $\varrho_{ 2-2\gamma }$ is defined by (\ref{lm 2.1}).

\end{lemma}
 {\it Proof.} Since problem (\ref{eq 1.1}) admits a super solution $u$ satisfying (\ref{2.3}) for  some $\lambda>0$, as in the proof of Proposition \ref{pr 2.1},  we know that
  (\ref{eq 1.1}) has the minimal solution $u_\lambda$ and then
$$   u_\lambda\le u\le \theta a,\qquad
  u_\lambda =\lambda \mathbb{G}_{\Omega}[\frac{1}{(a-u_\lambda)^2}]
   \le \lambda\kappa^{-2}(1-\theta)^{-2}\mathbb{G}_{\Omega}[\rho^{-2\gamma}],
$$
by the fact of  (\ref{1.1}).
Using  Lemma \ref{lm 2.1} with $\tau=2-2\gamma$, we   imply that (\ref{2.311}) holds.
The proof is complete. \qquad $\Box$

\begin{proposition}\label{pr 3.1.1}
Assume that the function  $a\in C^\gamma(\Omega)\cap C(\bar\Omega)$ satisfies (\ref{1.1}) and(\ref{a 1.0}) with $c_0\ge \kappa>0$ and $\gamma\in(0,\frac23]$.
Then \\
$(i)$  for $\lambda\in(0,\lambda^*)$, there exists $c_{21}\ge1$ dependent on $\gamma$ such that
$$\frac{\lambda}{c_{21}}\varrho_{ 2-2\gamma }(x)\le u_\lambda(x)\le c_{21}\lambda\rho(x)^{\gamma},\quad \forall x\in \Omega;$$
$(ii)$ there exists $\lambda_*\le \lambda^*$ such that for $\lambda\in(0,\lambda_*)$,
\begin{equation}\label{2.1.2}
 \frac{\lambda}{c_{21}}\varrho_{ 2-2\gamma }(x)\le u_\lambda(x)\le c_{21}\lambda\varrho_{ 2-2\gamma }(x),\quad \forall x\in \Omega
\end{equation}
 where $\varrho_{ 2-2\gamma }$ is defined by (\ref{lm 2.1}).
\end{proposition}
{\it Proof.} {\it Lower bound.}
By Proposition \ref{pr 2.1}, we see that  (\ref{eq 1.1}) admits the minimal solution $u_\lambda$ for $\lambda\in(0,\lambda^*)$,
which is approximated
by an increasing sequence $\{v_n\}_n$ defined by
$$ v_0=0 \quad{\rm and}\quad  v_n=\lambda\mathbb{G}_{\Omega}[\frac{1}{(a-v_{n-1})^2}].$$
 By the fact of (\ref{a 1.0}) and Lemma \ref{lm 2.1} with $\tau=2-2\gamma$, we have that
 \begin{eqnarray*}
 u_\lambda\ge v_1=\lambda\mathbb{G}_{\Omega}[\frac{1}{a^2}]
 \ge  \lambda c_0^2 \mathbb{G}_{\Omega}[\rho^{-2\gamma}]
    \ge   c_{22} \lambda\varrho_{ 2-2\gamma }\quad{\rm in}\quad \Omega
 \end{eqnarray*}
for some $c_{22}\in(0,1)$ dependent on $\gamma$.

{\it Upper bound.}  From the proof of Proposition \ref{pr 2.1}, for $\lambda>0$ small,  we see that $u_\lambda(x)\le \mu\rho(x)^\gamma$, $x\in\Omega$,
for some $\mu\in(0,\kappa)$, then there exists $\theta\in(0,1)$ such that
$$u_{\lambda}(x)\le \theta a(x),\quad x\in\Omega.$$
By Lemma \ref{lm 2.2}, we have that
\begin{equation}\label{eq921}
 u_\lambda\le c_{23}\varrho_{ 2-2\gamma }\quad{\rm in}\quad \Omega.
 \end{equation}
Let us define
$$\lambda_*=\sup\{\lambda\in(0,\lambda^*):\ \limsup_{x\in\Omega,x\to\partial\Omega}u_\lambda(x)\rho(x)^{-\gamma}<\kappa \},$$
we observe that $\lambda_*\le \lambda^*$ and  (\ref{eq921}) holds for all $\lambda\in(0,\lambda_*)$.
This ends the proof.
\qquad $\Box$

\setcounter{equation}{0}
\section{Estimates for $\lambda^*$ and $\lambda_*$ when $\Omega=B_1(0)$ }

In this section, we do the estimates for $\lambda^*$ and $\lambda_*$ in the case that $\Omega=B_1(0)$
and
\begin{equation}\label{a 1.1}
a(x)=\kappa(1-|x|^2)^{\gamma}.
\end{equation}
In this case, the function $a$ represents the upper semi-sphere type shape in $\R^N$. We have following monotonicity results.

\begin{proposition}\label{lm 3.1}
Assume that $\Omega=B_1(0)$, the function
$a$ satisfies (\ref{a 1.1})
with $\kappa>0$, $\gamma\in(0,\frac23]$ and $0<\lambda<\lambda^*(\kappa,\gamma)$. Then

$(i)$ the mappings: $\gamma\mapsto \lambda^*(\kappa,\gamma)$ and $\gamma\mapsto \lambda_*(\kappa,\gamma)$  are decreasing for fixed  $\kappa>0$;

$(ii)$ the mapping: $\kappa\mapsto \lambda^*(\kappa,\gamma)$  and $\kappa\mapsto \lambda_*(\kappa,\gamma)$ are increasing for fixed $\gamma\in(0,\frac23]$.
\end{proposition}
{\it Proof.}  Let $0<\gamma_2\le \gamma_1\le \frac23$, for $\lambda\in(0,\min\{\lambda^*(\kappa, \gamma_1), \lambda^*(\kappa, \gamma_2)\})$,
 denote by $u_1, u_2$  the minimal solutions of (\ref{eq 1.1})
with $a=a_1$ and $a=a_2$ respectively, where $a_1(x)=\kappa (1-|x|^2)^{\gamma_1}$ and $a_2(x)=\kappa (1-|x|^2)^{\gamma_2}$. Since $\gamma_2\le \gamma_1$,
then $a_1\le a_2$ in $B_1(0)$
and  for any $\lambda\in(0,\lambda^*(\kappa, \gamma_1))$,
$$-\Delta u_1=\frac{\lambda}{(a_1-u_1)^2}\ge \frac{\lambda}{(a_2-u_1)^2}\quad {\rm in} \ \ B_1(0),$$
that is, $u_1$ is a super solution of (\ref{eq 1.1})
with  $a=a_2$, for any  $\lambda\in(0,\lambda^*(\kappa, \gamma_1))$. By the same argument in Proposition \ref{pr 2.1},
problem (\ref{eq 1.1})
with  $a=a_2$ admits the minimal solution for any $\lambda\in(0,\lambda^*(\kappa, \gamma_1))$. Thus,
$$\lambda^*(\kappa, \gamma_2)\ge \lambda^*(\kappa, \gamma_1),$$
by the definition of $\lambda^*(\kappa, \gamma_2))$. As a consequence, we obtain that
the mapping $\gamma\mapsto \lambda^*(\kappa,\gamma)$ is decreasing for fixed  $\kappa>0$.

It is similar to obtain the other assertions. \qquad$\Box$

\begin{proposition}\label{pr 3.2}
Assume that $\Omega=B_1(0)$ and  the function
$a$ satisfies (\ref{a 1.1})
with $\kappa>0$, $\gamma\in(0,\frac23]$.
 Then
\begin{equation}\label{3.1-1}
 \lambda^*(\kappa,\gamma)\le \frac{4N\kappa^3}{3}\mathcal{B}(\frac12,2-2\gamma),
\end{equation}
where $\mathcal{B}(\cdot,\cdot)$ is the Beta function.
 \end{proposition}
{\it Proof.} From (\ref{2.4}), we  only need to do the estimates of $\int_{B_1(0)}a(x)dx$ and $ \int_{B_1(0)}\frac{\mathbb{G}_\Omega[1](x)}{a^2(x)}dx$.
Since $a$ satisfies (\ref{a 1.1}),  by direct computation, we have that
\begin{equation}\label{a b1}
\int_{B_1(0)}a(x)dx=  \frac23\kappa|S^{N-1}|,
\end{equation}
where $S^{N-1}$ is the unit sphere in $\R^N$.

On the other hand, we see that  there exists $c_{27}>1$ such that
 $$ \mathbb{G}_{B_1(0)}[1](x)=  \frac1{2N}(1-|x|^2),\quad \forall x\in B_1(0),$$
 then
 \begin{eqnarray*}
 \int_{B_1(0)}\frac{\mathbb{G}_\Omega[1](x)}{a^2(x)}dx= \frac1{2N\kappa^2} \int_{B_1(0)}(1-|x|^2)^{1-2\gamma} dx  =  \frac{|S^{N-1}|}{4N\kappa^2} \mathcal{B}(\frac12,2-2\gamma),
 \end{eqnarray*}
which, combining with (\ref{a b1}),  (\ref{2.4}), implies (\ref{3.1-1}).
\qquad$\Box$

\smallskip

From Proposition \ref{lm 3.1}, the mapping $\gamma\mapsto \lambda_*(\kappa,\gamma)$  is decreasing,
so our interest is to evaluate   the lower bound for $\lambda_*$ when $\gamma=\frac{2}{3}$.

\begin{proposition}\label{pr 3.1}
Assume that $\Omega=B_1(0)$ and the function $a$ satisfies (\ref{a 1.1})
with $\kappa>0$ and $\gamma=\frac23$.
Then
\begin{equation}\label{3.3}
 \lambda_*(\kappa,\frac23)\ge
\left(\frac{2}{3}\right)^5\kappa^3.
\end{equation}

\end{proposition}
{\it Proof.} Let $w_t(r)=t\kappa (1-r^2)^{2/3}$ with $t\in(0,1)$, then for $x\in B_1(0)$,
\begin{eqnarray*}
-\Delta w_t(|x|)&=&  \frac89 t\kappa  (1-|x|^2)^{-4/3}+ \frac43t\kappa   (N-\frac23) (1-|x|^2)^{ -1/3} \\
   &\ge& \frac89 t\kappa  (1-|x|^2)^{-4/3}
\end{eqnarray*}
and
$$
\frac{\lambda }{(a(x)-w(|x|))^2}\le (1-t)^{-2}\kappa^{-2} \lambda   (1-|x|^2)^{-4/3},
$$
thus, if
\begin{equation}\label{3.1-2}
 \lambda  \le \lambda_t:=\frac89 t(1-t)^2\kappa^3,
\end{equation}
$w_t$ is a super solution of the equation (\ref{eq 1.1}).
In fact, choose $t=\frac13$, we get the optimal $\lambda_t=\frac{2^5}{3^4}\kappa^3$ and then by
 Lemma \ref{lm 2.2} to obtain that problem (\ref{eq 1.1}) admits the minimal solution $u_{\lambda_t}$. Therefore,
 $\lambda_*\ge \lambda_t$. The proof is complete.
\qquad$\Box$

\smallskip

\begin{lemma}\label{lm 4.4}
Assume that
$a\in C^\gamma(\Omega)\cap C(\bar\Omega)$  satisfies  (\ref{1.1}) and (\ref{a 1.0})
with $c_0\ge\kappa>0$, $\gamma\in(0,\frac23]$, $\lambda_*$ is given in Theorem \ref{teo 1} and $ u_\lambda$ is the minimal solution of (\ref{eq 1.1})
with  $\lambda\in(0,\lambda_*)$.
Then $u_\lambda\in H^1_0(\Omega)$ and  there exists $c_{31}>0$   such that
$$ \int_{\Omega} \frac{u_\lambda}{(a-u_{\lambda})^2} dx\le c_{31}\lambda \qquad{\rm and}\qquad \int_\Omega |\nabla u_\lambda|^2 dx\le c_{31}\lambda^2.$$

\end{lemma}
{\it Proof.} For $\lambda\in(0,\lambda_*)$ and $\gamma\in(0,\frac23]\setminus{\{\frac12\}}$, by (\ref{2.1.2}), we have that
$$ u_\lambda(x)\le c_{21}\lambda\rho(x)^{\min\{1,2-2\gamma\}},\quad \  x\in\Omega.$$
Since $\kappa \rho(x)^\gamma\le a(x)\le c_0\rho(x)^\gamma$ for $x\in\Omega$,
then there exists $\theta_1\in(0,1)$ such that $u_\lambda<\theta_1 a$ in $A_{\delta}=\{x\in\Omega:\, \rho(x)<\delta\}$ for some $\delta>0$ small.
Moreover, since $u_\lambda<a $ in $\Omega$, then there exists $\theta_2\in(0,1)$ such that $u_\lambda\le \theta_2a$ in $\Omega\setminus{A_{\delta}}$.
Let $\theta=\max\{\theta_1,\theta_2\}$, then we have that $u_\lambda\le \theta a$ in $\Omega$ and then
\begin{equation}\label{eq 01}
a-u_\lambda\ge (1-\theta) a \ge (1-\theta)\kappa \rho^\gamma\quad {\rm in} \ \ \Omega.
 \end{equation}
Therefore,
 \begin{equation}\label{4.1.2}
 \int_\Omega \frac{u_\lambda}{(a-u_\lambda)^2}dx\le \int_\Omega \frac{ c_{21}\lambda\rho(x)^{\min\{1,2-2\gamma\}}}{(1-\theta)^2\kappa^2 \rho(x)^{2\gamma}}dx:=c_{24}\lambda<+\infty,
 \end{equation}
 by the fact that $\min\{1,2-2\gamma\}-2\gamma>-1$. Taking a sequence $\{\xi_n\}_n\subset C_c^2(\Omega)$  which converges to
$u_\lambda$ as $n\to\infty$, since $u_\lambda$ is the minimal solution of  (\ref{eq 1.1}), then
$$\int_\Omega \nabla u_\lambda\cdot \nabla{\xi_n} dx=\int_\Omega \frac{\lambda \xi_n}{(a-u_\lambda)^2} dx.$$
Passing to the limit as $n\to\infty$, we obtain that
$$
 \int_\Omega  |\nabla u_\lambda|^2dx =  \int_\Omega \frac{\lambda u_\lambda}{(a-u_\lambda)^2}dx\le c_{24}\lambda^2.
$$

For $\lambda\in(0,\lambda_*)$ and $\gamma=\frac12$,   by (\ref{2.1.2}), we have that
$$ u_\lambda(x)\le c_{21}\lambda\rho(x)\log\frac{1}{\rho(x)},\quad \  x\in\Omega,$$
similarly, there exists $\theta\in(0,1)$ such that $u_\lambda\le \theta a$ in $\Omega$ and then (\ref{eq 01}) holds. Then
 \begin{equation}\label{4.1.20p}
 \int_\Omega \frac{u_\lambda}{(a-u_\lambda)^2}dx\le \int_\Omega \frac{ c_{21}\lambda}{(1-\theta)^2\kappa^2 }\log\frac{1}{\rho(x)}dx\le c_{24}\lambda
 \end{equation}
and then we have that $ \int_\Omega  |\nabla u_\lambda|^2dx \le c_{24}\lambda^2.$
The proof is complete.\qquad$\Box$

\medskip

\noindent{\bf Proof of Theorem \ref{teo 1}.} The existence of the minimal solution for $\lambda\in(0,\lambda^*)$
and the nonexistence for $\lambda>\lambda^*$ follow by Proposition \ref{pr 2.1}. The proof of  Theorem \ref{teo 1} $(iii)$
sees Proposition \ref{pr 3.1.1} and Lemma \ref{lm 4.4}. The estimates of $\lambda^*$ and $\lambda_*$ see Proposition \ref{lm 3.1},
Proposition \ref{pr 3.1} and Proposition \ref{pr 3.2}.\qquad $\Box$

\setcounter{equation}{0}
\section{Properties of minimal solution}

\subsection{Regularity}

In this subsection, we study the solutions of (\ref{eq 1.1}) when $\lambda=\lambda^*$.
\begin{proposition}\label{pr 4.1}

Assume that the function
$a\in C^\gamma(\Omega)\cap C(\bar\Omega)$  satisfies (\ref{1.1}) and (\ref{a 1.0})
with $c_0\ge\kappa>0$, $\gamma\in(0,\frac23]$ and $ u_{\lambda^*}$ is given by (\ref{4.2}).
Then  $ u_{\lambda^*}$ is a weak solution of (\ref{eq 1.1}) with $\lambda=\lambda^*$. Moreover,
for any $\beta\in(0,\gamma)$,  there exists $c_\beta>0$ such that
 \begin{equation}\label{4.0.5}
 \norm{u_{\lambda^*}}_{W^{1,\frac{N}{N-\beta}}(\Omega)}\le c_\beta
 \end{equation}
and
\begin{equation}\label{4.0.6}
 \int_{\Omega} \frac{\rho^{1-\beta}}{(a-u_{\lambda^*})^2} dx\le c_\beta.
 \end{equation}
\end{proposition}
{\it Proof.} For any $\beta\in(0,\gamma)$ and $n\in\N$, denote $\xi_n=\mathbb{G}_{\Omega}[\rho^{-1-\beta}]\eta_n$, where $\eta_n$ is defined by (\ref{101}),
 we observe that $\xi_n\in C_c^2(\Omega)$ and
\begin{eqnarray*}
\int_\Omega u_\lambda (-\Delta\xi_n) dx=\int_\Omega u_\lambda ( \rho^{-1-\beta}\eta_n-\nabla \mathbb{G}_{\Omega}[\rho^{-1-\beta}]\cdot\nabla\eta_n+(-\Delta)\eta_n\mathbb{G}_{\Omega}[\rho^{-1-\beta}])dx.
\end{eqnarray*}
Using Lemma \ref{lm 2.1} with $\tau=1-\beta\in(1/3,1)$, we have that $\mathbb{G}_{\Omega}[\rho^{-1-\beta}]\le c_{25}\rho^{1-\beta}$ in $\Omega$.
Combining with the fact that $|\nabla \eta_n|\le c_{16}n$, $|(-\Delta)\eta_n|\le c_{17}n^2$
and  $0<u_\lambda(x)<a(x)\le c_0\rho(x)^\gamma$ for $x\in\Omega$, we obtain that
\begin{eqnarray*}
  \int_\Omega u_\lambda (-\Delta\xi_n) dx&\le& c_0\int_\Omega  \rho^{\gamma-1-\beta} dx+c_{26}n\int_{A_{1/n}}  \rho^{\gamma-\beta}dx + c_{26}n^2\int_{A_{1/n}} \rho^{1-\beta+\gamma} dx \\
   &\le & \bar c_{\beta},
\end{eqnarray*}
where $ c_{26}, \bar c_{\beta}>0$ independent on $n$ and $\bar c_\beta$ satisfies $\bar c_\beta\to+\infty$ as $\beta\to\gamma^-$. Thus,
$$\int_\Omega \frac{ \lambda  \mathbb{G}_{\Omega}[\rho^{-1-\beta}]}{(a-u_\lambda)^2} dx\le \bar c_\beta$$
and then
\begin{equation}\label{4.6p}
\int_\Omega \frac{ \lambda  \rho^{1-\beta}}{(a-u_\lambda)^2} dx\le \frac{\bar c_\beta}{c_{25}}.
\end{equation}
By Theorem \ref{teo 1}, we see that the mapping $\lambda\mapsto u_\lambda$ is increasing and
uniformly bounded by the function $a$, which is in $L^1(\Omega)$, then
$$u_\lambda\to u_{\lambda^*}\ \ \  {\rm in}\ \ L^1(\Omega)\quad \ {\rm as}\ \lambda\to \lambda^*$$
and then for any $\xi\in C_c^2(\Omega)$, we have that
\begin{equation}\label{4.6p1}
\int_\Omega u_\lambda(-\Delta\xi) dx\to \int_\Omega u_{ \lambda^*}(-\Delta \xi) dx \ \ \ {\rm as} \ \ \lambda\to  \lambda^*.
\end{equation}
Moreover,  we observe that the mapping $\lambda\mapsto  \frac{\lambda  }{(a-u_\lambda)^2}$ is increasing and by (\ref{4.6p}),
$$\frac{\lambda  }{(a-u_\lambda)^2}\to \frac{\lambda^*  }{(a-u_{\lambda^*})^2}\ \ \  {\rm a.e.\ in }\ \  \Omega\quad {\rm as}\ \lambda\to \lambda^*,$$
Then  (\ref{4.6p}) deduces that
$$\frac{\lambda  }{(a-u_\lambda)^2}\to \frac{\lambda^*  }{(a-u_{\lambda^*})^2}\ \ \  {\rm in}\ \ L^1(\Omega,\,\rho^{1-\beta}dx)\quad {\rm as}\ \lambda\to \lambda^*$$
and
\begin{equation}\label{4.6p3}
\int_\Omega \frac{ \lambda^*  \rho^{1-\beta}}{(a-u_\lambda^*)^2} dx\le \frac{\bar c_\beta}{c_{25}},
\end{equation}
thus, for any $\xi\in C_c^2(\Omega)$, we have that
\begin{equation}\label{4.6p2}
\int_\Omega \frac{\lambda \xi }{(a-u_\lambda)^2} dx\to \int_\Omega \frac{\lambda^* \xi }{(a-u_{\lambda^*})^2} dx \ \ \ {\rm as} \ \ \lambda\to  \lambda^*.
\end{equation}
Since $u_\lambda$ is the minimal solution of  (\ref{eq 1.1}),
$$
\int_\Omega u_\lambda(-\Delta\xi) dx=\int_\Omega \frac{\lambda \xi}{(a-u_\lambda)^2} dx,\ \ \ \ \ \forall \, \xi\in C_c^2(\Omega),
$$
 passing to the limit as $\lambda\to\lambda^*$, combining with (\ref{4.6p1}) and (\ref{4.6p2}),
then
 $u_{\lambda^*}$ is a weak solution of (\ref{eq 1.1}) with $\lambda=\lambda^*$.

By \cite[Theorem 2.6]{BV},   for any $\beta'\in(\beta,\gamma)$,  there exists $c_{27}>0$ such that
$$
\norm{|\nabla u_{\lambda^*}|}_{M^{\frac{N}{N- \beta'}}(\Omega) }\le c_{27}\norm{(a-u_{\lambda^*})^{-2}}_{L^1(\Omega, \, \rho^{1-\beta'}dx)},
$$
where ${M^{\frac{N}{N- \beta'}}(\Omega) }$ is the Marcinkiewicz space of
exponent $\frac{N}{N- \beta'}$.
Then, by (\ref{4.6p3}), we have that
$$\norm{|\nabla u_{\lambda^*}|}_{L^{\frac{N}{N- \beta}}(\Omega)}\le c_{27}c_{25}^{-1}\bar c_\beta(\lambda^*)^{-1}$$
and then
$$
 \norm{u_\lambda^*}_{W^{1,\frac{N}{N- \beta}}(\Omega)}\le c_\beta
$$
for some $c_\beta>0.$
Thus, (\ref{4.0.5}) holds and  (\ref{4.0.6}) is obvious by (\ref{4.6p3}).
This ends the proof. \qquad$\Box$

\subsection{Stability}

In this subsection, we introduce the stability of the minimal solution $u_\lambda$ for problem (\ref{eq 1.1}). By the definition of $\lambda_*$, we observe that for
any $\lambda\in (0,\lambda_*)$, there exists $\theta\in(0,1)$ such that $u_\lambda\le \theta  a$. Therefore,
by the fact that $\gamma\in(0,\frac23]$ and $a\ge \kappa \rho^\gamma$ in $\Omega$, we have that
\begin{equation}\label{5.1}
\frac{1}{(a-u_\lambda)^3}\le \frac{1}{(1-\theta)^3a^3}\le (1-\theta)^{-3}\kappa^{-3}\rho^{-3\gamma}\le   c_{28}\rho^{-2}\quad{\rm in}\quad \Omega,
\end{equation}
where $c_{28}>0$ depends on $\theta, \kappa, \gamma$.
This enables us to consider the first eigenvalue of $-\Delta-\frac{2\lambda}{(a-u_\lambda)^3}$ in $H_0^1(\Omega)$, that is,
$$\mu_1(\lambda)=\inf_{\varphi\in H_0^1(\Omega)}\frac{\int_\Omega (|\nabla \varphi|^2-\frac{2\lambda \varphi^2}{(a-u_\lambda)^3})dx}{\int_\Omega \varphi^2\,dx}$$
for $\lambda\in(0,\lambda_*)$. It is well-known that  $u_\lambda$ is stable if $\mu_1(\lambda)>0$  and semi-stable if $\mu_1(\lambda)\ge0$.

\begin{lemma}\label{lm 4.1}
Assume that $\lambda\in(0,\lambda_*)$,
$a\in C^\gamma(\Omega)\cap C(\bar\Omega)$ satisfies (\ref{1.1}) and (\ref{a 1.0})
with $c_0\ge\kappa>0$, $\gamma\in(0,\frac23]$. Suppose that $u_\lambda$ is the minimal  solution of (\ref{eq 1.1})
and $v_\lambda$ is a super solution of (\ref{eq 1.1}).

If $\mu_1(\lambda)>0$, then
$$u_\lambda\le v_\lambda\quad{\rm in}\quad \Omega.$$

If $\mu_1(\lambda)=0$, then
$$u_\lambda= v_\lambda\quad{\rm in}\quad \Omega.$$
\end{lemma}
{\it Proof.}  We observe that $a>u_\lambda$ in $\Omega$ for $\lambda\in(0,\lambda^*)$ and $\mu_1(\lambda)$ is well-defined for $\lambda\in(0,\lambda_*)$, then it follows the procedure of the proof of Lemma 4.1 in \cite{GG} just replacing
$\frac{f}{(1-u)^2}$ by $\frac{1}{(a-u)^2}$.

\hfill$\Box$

\begin{proposition}\label{pr 4.2}
Assume that $\lambda\in(0,\lambda_*)$,
$a\in C^\gamma(\Omega)\cap C(\bar\Omega)$ satisfies  (\ref{1.1}) and (\ref{a 1.0})
with $c_0\ge\kappa>0$, $\gamma\in(0,\frac23]$. Let $u_\lambda$ be the minimal solution of (\ref{eq 1.1}),
then $u_\lambda$ is stable.

\end{proposition}
{\it Proof.}
  We first claim  that the mapping $\lambda\mapsto \mu_1(\lambda)$ is locally Lipschitz continuous and strictly decreasing
in $(0,\lambda_*)$.
Observing that the mapping $\lambda\mapsto u_{\lambda}$ is strictly increasing and so is $\lambda\mapsto \frac{2\lambda}{(a-u_\lambda)^3}$,
which implies that the mapping $\lambda\mapsto \mu_1(\lambda)$ is strictly decreasing
in $(0,\lambda_*)$. Let $0<\lambda_1<\lambda_2<\lambda_*$ and $\phi_{\lambda_2}$ be the achieved function of $\mu_1(\lambda_2)$ with $\norm{\phi_{\lambda_2}}_{L^2(\Omega)}=1$,
then we have that
\begin{eqnarray*}
0< \mu_1(\lambda_1)-\mu_1(\lambda_2) &\le& \int_\Omega \left(|\nabla \phi_{\lambda_2}|^2-\frac{2\lambda_1 \phi_{\lambda_2}^2}{(a-u_{\lambda_1})^3}\right)dx  \\
   && -\int_\Omega \left(|\nabla \phi_{\lambda_2}|^2-\frac{2\lambda_2 \phi_{\lambda_2}^2}{(a-u_{\lambda_2})^3}\right)dx
   \\&\le &2(\lambda_2-\lambda_1)\int_\Omega\frac{ \phi_{\lambda_2}^2}{(a-u_{\lambda_2})^3} dx,
\end{eqnarray*}
thus, the mapping $\lambda\mapsto \mu_1(\lambda)$ is locally Lipschitz continuous.

 We next prove that $u_\lambda$ is stable and $\mu_1(\lambda)>0$ for $\lambda>0$ small.
It follows by \cite[Theorem 1]{MS} that there exists constant $c_{29}>0$ such that
\begin{equation}\label{4.1.1}
\int_\Omega \varphi^2\rho^{-2}dx\le c_{29}\int_{\Omega} |\nabla \varphi|^2 dx,\quad \forall \varphi\in H^1_0(\Omega).
\end{equation}
For $\lambda<\lambda_*$, it follows from
 (\ref{5.1}) and (\ref{4.1.1}) that
$$\int_{\Omega}\frac{2\lambda \varphi^2}{(a-u_\lambda)^3}dx\le 2\lambda c_{28}\int_\Omega \varphi^2\rho^{-2}dx
\le 2\lambda c_{28}c_{29} \int_{\Omega} |\nabla \varphi|^2 dx, \quad \forall \varphi\in H^1_0(\Omega).$$
For   $\lambda>0$ small such that $2\lambda c_{28}c_{29}< 1$, we obtain that
$$\int_{\Omega}\frac{2\lambda \varphi^2}{(a-u_\lambda)^3}dx
<  \int_{\Omega} |\nabla \varphi|^2 dx, \quad \forall \varphi\in H_0^1(\Omega)\setminus\{0\},$$
that is, $\mu_1(\lambda)>0$ and $u_\lambda$ is  stable   if $\lambda>0$ small.

 Finally, we prove that  $u_{\lambda}$ is stable for $\lambda\in(0,\lambda_*)$.    We have obtained that $\mu_1(\lambda)>0$ for $\lambda>0$ small and
 the mapping $\lambda\mapsto \mu_1(\lambda)$ is locally Lipschitz continuous, strictly decreasing
in $(0,\lambda_*)$, so if there exists $\lambda_0\in(0,\lambda_*)$ such that $\mu_1(\lambda_0)=0$, then $\mu_1(\lambda)>0$ for $\lambda\in(0,\lambda_0)$.
Letting $\lambda_1\in (\lambda_0,\lambda_*)$, the minimal solution $u_{\lambda_1}$  is
a super solution of
$$
\left\{\arraycolsep=1pt
\begin{array}{lll}
 -\Delta    u = \frac{\lambda_0 }{(a-u)^2}\quad  &{\rm in}\quad\ \Omega,
 \\[2mm]
 \phantom{-\Delta   }
 u=0\quad &{\rm on}\quad   \partial \Omega
\end{array}
\right.
$$
and it infers from Lemma \ref{lm 4.1} that
$$u_{\lambda_1}=u_{\lambda_0},$$
 which contradicts that the mapping $\lambda\mapsto u_{\lambda}$ is strictly increasing.
 Therefore, $\mu_1(\lambda)>0$  and  $u_{\lambda}$ is stable for $\lambda\in(0,\lambda_*)$.
 \qquad$\Box$

\subsection{Particular case that $\gamma=\frac23$}

In this subsection, we study further on properties of the minimal solution  when $\gamma=\frac23$.

\begin{lemma}\label{lm 5.1}
Assume that
$a\in C^\gamma(\Omega)\cap C(\bar\Omega)$ satisfies  (\ref{1.1}) and (\ref{a 1.0})
with $c_0\ge\kappa>0$, $\gamma=\frac23$. Then $\lambda^*=\lambda_*$.

\end{lemma}
{\it Proof.} For $\lambda\in(0,\lambda^*)$,
by (\ref{e 1.1}) and Lemma \ref{lm 2.1}, it implies that
 \begin{eqnarray*}
  a-u_{\lambda}& \ge& (\lambda^*-\lambda)\mathbb{G}_\Omega[a^{-2}]\\&\ge& c_{30}(\lambda^*-\lambda)\rho^{\gamma}
    \ge  c_{31}(\lambda^*-\lambda)a,
 \end{eqnarray*}
 then there exists $\theta\in(0,1)$ such that $u_\lambda\le \theta a$ in $\Omega$. It follows by Lemma \ref{lm 2.2} and the definition of $\lambda_*$ that $\lambda_*=\lambda^*$.
\qquad$\Box$

\begin{proposition}\label{pr 4.3}

Assume that
$a\in C^\gamma(\Omega)\cap C(\bar\Omega)$ satisfies  (\ref{1.1}) and (\ref{a 1.0})
with $c_0\ge\kappa>0$,  $\gamma=\frac23$ and $ u_{\lambda^*}$ is given by (\ref{4.2}).
Then  $ u_{\lambda^*}$ is a semi-stable weak solution of (\ref{eq 1.1}) with $\lambda=\lambda^*$.

\end{proposition}
{\it Proof.} When  $\gamma=\frac23$,
by Lemma \ref{lm 5.1} and Proposition \ref{pr 4.2}, we have that  $u_\lambda$ is stable for $\lambda\in(0,\lambda^*)$, that is,
$$
\int_\Omega\frac{2\lambda\varphi^2}{(a-u_\lambda)^3}dx<\int_\Omega|\nabla \varphi|^2dx,\quad \ \forall \varphi\in H^1_0(\Omega)\setminus\{0\}.
$$
Let $\varphi=\mathbb{G}_{\Omega}[1]$, we have that
$$\int_\Omega\frac{ \rho^2}{(a-u_\lambda)^3}dx<c_{39}\lambda^{-1}.$$
Since the mapping $\lambda\mapsto \frac{ \rho^2}{(a-u_\lambda)^3}$ is strictly increasing
and bounded in $L^1(\Omega)$,
then
$$\frac{ \rho^2}{(a-u_\lambda)^3}\to \frac{ \rho^2}{(a-u_{\lambda^*})^3}\quad{\rm as}\quad \lambda\to\lambda^*\quad{\rm in}\quad L^1(\Omega)$$
and for any $\varphi\in C_c^2(\Omega)$,
$$\lim_{\lambda\to \lambda^*}\int_\Omega\frac{\lambda\varphi^2}{(a-u_\lambda)^3}dx=\int_\Omega\frac{\lambda^*\varphi^2}{(a-u_{\lambda^*})^3}dx$$
Therefore,
$$\int_\Omega\frac{2\lambda^*\varphi^2}{(a-u_{\lambda^*})^3}dx\le \int_\Omega|\nabla \varphi|^2dx,\quad \forall \varphi\in C_c^2(\Omega),$$
by the fact that $C_c^2(\Omega)$ is dense in $H^1_0(\Omega)$,
then $u_{\lambda^*}$ is semi-stable.
\qquad$\Box$

\medskip

We next improve the regularity of $u_{\lambda^*}$ and prove when $1\le N\le 7$, the extremal
solution $u_{\lambda^*}$ is a classical solution of (\ref{eq 1.1}) with $\lambda=\lambda^*$. To this end,
we need the following Lemma, which is inspired by \cite{GG}.

\begin{lemma}\label{lm 4.2}
Assume that $\lambda\in(0,\lambda^*)$, $\Omega=B_1(0)$ and $a(x)=\kappa(1-|x|^2)^{\frac23}$
with $\kappa>0$.
Let $u$ be a weak solution of (\ref{eq 1.1}) satisfying for any compact set $K\subset B_1(0)$, there exists
$c_{32}>0$ such that
 \begin{equation}\label{4.1.3}
 \norm{\frac1{a-u}}_{L^{\frac{3N}2}(K)}\le c_{32}.
 \end{equation}
Then there exists $c_{33}>0$ depending on $K$ such that
\begin{equation}\label{4.2.1}
 \inf_{x\in K}(a(x)-u(x))>c_{33}.
\end{equation}
\end{lemma}
{\it Proof.}  By (\ref{4.1.3}), we have that
$$\frac1{(a-u)^2}\in L^{\frac{3N}4}(K)$$
and then $u\in W^{2,\frac{3N}4}(K)$ and by Sobolev's Theorem, we
deduce that $u\in C^{\frac23}(K')$ with $K'$  compact set in interior point set of $K$. We next
show that $u<a$ in $B_1(0)$. Indeed, if not, there exists $x_0\in B_1(0)$ such that $u(x_0)=a(x_0)$. For compact set $K\subset B_1(0)$ containing $x_0$,
  we have that
\begin{eqnarray*}
 |a(x)-u(x)|&\le &|a(x)-a(x_0)|+|u(x)-u(x_0)|+|u(x_0)-a(x_0)|\\ &=&
|a(x)-a(x_0)|+|u(x)-u(x_0)|\le c_{34}|x-x_0|^{\frac23},
\end{eqnarray*}
then
$$\int_{K}\frac1{(a-u)^{\frac{3N}2}}dx\ge c_{35} \int_{K} |x-x_0|^{-\frac{3N}2\cdot\frac23}dx=\infty,$$
which  contradicts (\ref{4.1.3}). Therefore, (\ref{4.2.1}) holds.
 \qquad$\Box$

\begin{proposition}\label{pr 4.4}

Assume that $1\le N\le 7$, $\Omega=B_1(0)$,
 the function
$a$ satisfies (\ref{a 1.1})
with $\kappa>0$ and $\gamma=\frac23$, $ u_{\lambda^*}$ is given by (\ref{4.2}).
Then  $u_{\lambda^*}$ is a classical solution of (\ref{eq 1.1}) with $\lambda=\lambda^*$.

\end{proposition}
{\it Proof.}  Since the mapping $\lambda\mapsto u_\lambda$ is strictly increasing and bounded by $a$, then from (\ref{a 1.0}) and Lemma \ref{lm 4.2},
we only have to improve the regularity of $ u_{\lambda^*}$ in any compact set of $B_1(0)$. For $\lambda\in(0,\lambda^*)$, we know that $u_\lambda$ is stable, then
\begin{equation}\label{4.3.1}
\int_{B_1(0)}\frac{2\lambda \varphi^2}{(a-u_\lambda)^3}dx< \int_{\Omega} |\nabla \varphi|^2 dx, \quad \forall \varphi\in H^1_0(B_1(0))\setminus\{0\}.
\end{equation}

We claim that the minimal solutions $u_{\lambda}$ is radially symmetric for $\lambda\in(0,\lambda^*]$. Indeed,    the minimal solution $u_\lambda$ could be approximated
by the sequence of functions
$$v_n=\lambda\mathbb{G}_{B_1(0)}[\frac1{(a-v_{n-1})^2}]\ \ \ {\rm with} \ \  v_0=0.$$
It follows by radially symmetry of $v_{n-1}$ and $a$ that $v_n$ is radially symmetry and then
$u_\lambda$ is radially symmetric. Then $u_{\lambda^*}$ is radially symmetric.

By (\ref{4.0.6}),   there exists a  sequence $\{r_n\}_n\subset(0,1)$ such that
$$\lim_{n\to+\infty} r_n=1\quad{\rm and}\quad  a(r_n)-u_{\lambda^*}(r_n)>0.$$
Otherwise, there exists $r\in(0,1)$ such that
$$a-u_{\lambda^*}=0\quad{\rm a.e.\ in}\ \ B_1(0)\setminus{B_r(0)},$$
which contradicts (\ref{4.0.6}).

Let us denote
$$
\varphi_\theta= \left\{\arraycolsep=1pt
\begin{array}{lll}
(a-u_\lambda)^\theta-\epsilon_\lambda^\theta\quad  &{\rm in}\quad B_{r_n}(0),
 \\[1.5mm]
 \phantom{}
0\quad  &{\rm in}\quad B_1(0)\setminus B_{r_n}(0),
\end{array}
\right.
$$
where $\theta\in(-2-\sqrt{6},0)$ and $\epsilon_\lambda= a(r_n)-u_\lambda(r_n)$. We observe that $\varphi_\theta\in H_0^1(B_1(0))$.
It follows by (\ref{4.3.1}) with $\varphi_\theta$, we have that
\begin{eqnarray}
  \int_{B_{r_n}(0)}\frac{2\lambda [(a-u_\lambda)^\theta-\epsilon_\lambda^\theta]^2}{(a-u_\lambda)^{3}}dx &\le &
 \int_{B_{r_n}(0)}|\nabla ((a-u_\lambda)^\theta)|^2dx\nonumber \\
   &=& \theta^2\int_{B_{r_n}(0)}(a-u_\lambda)^{2\theta-2}  |\nabla (a-u_\lambda)|^2  dx.\label{4.3.2}
\end{eqnarray}
Since $u_\lambda$ is the minimal solution of (\ref{eq 1.1}), then
\begin{equation}\label{4.3.3}
  -\Delta (a-u_\lambda)=-\Delta a-\frac{\lambda}{(a-u_\lambda)^2}\quad {\rm in} \ \  B_{r_n}(0).
\end{equation}
Multiplying  (\ref{4.3.3})  by $\frac{\theta^2}{1-2\theta} [(a-u_\lambda)^{2\theta-1}-\epsilon_\lambda^{2\theta-1}]$ and applying integration by parts yields that
\begin{eqnarray*}
&&\frac{\theta^2}{1-2\theta}\int_{B_{r_n}(0)}[-\Delta a-\frac{\lambda}{(a-u_\lambda)^2}][(a-u_\lambda)^{2\theta-1}-\epsilon_\lambda^{2\theta-1}]dx
\\&&\qquad=\frac{\theta^2}{1-2\theta}\int_{B_{r_n}(0)}\nabla (a-u_\lambda)\cdot\nabla\left((a-u_\lambda)^{2\theta-1}\right)dx
 \\  &&\qquad\qquad =  -\theta^2 \int_{B_{r_n}(0)}(a-u_\lambda)^{2\theta-2}|\nabla(a-u_\lambda)|^2dx,
\end{eqnarray*}
 then together with (\ref{4.3.2}), we deduce that
\begin{eqnarray*}
\int_{B_{r_n}(0)}\frac{2\lambda [(a-u_\lambda)^\theta-\epsilon_\lambda^\theta]^2}{(a-u_\lambda)^{3}}dx \le   \frac{\theta^2}{1-2\theta}\int_{B_{r_n}(0)}[ \Delta a+\frac{\lambda}{(a-u_\lambda)^2}][(a-u_\lambda)^{2\theta-1}-\epsilon_\lambda^{2\theta-1}]dx,
\end{eqnarray*}
thus,
\begin{eqnarray*}
\lambda(2-\frac{\theta^2}{1-2\theta})\int_{B_{r_n}(0)}\frac1{(a-u_\lambda)^{3-2\theta}}dx\le \int_{B_{r_n}(0)}\frac{4\lambda\epsilon_\lambda^\theta}{(a-u_\lambda)^{3-\theta}}dx
-\int_{B_{r_n}(0)}\frac{2\lambda\epsilon_\lambda^{2\theta}}{(a-u_\lambda)^{3}}dx
\\ \ \  +\frac{\theta^2}{1-2\theta}\int_{B_{r_n}(0)} \frac{\Delta a}{(a-u_\lambda)^{1-2\theta}}dx-\frac{\theta^2\epsilon_\lambda^{2\theta-1}}{1-2\theta}\int_{B_{r_n}(0)}\Delta a dx
  -\frac{\theta^2\lambda}{1-2\theta}\int_{B_{r_n}(0)}\frac{\epsilon_\lambda^{2\theta-1}}{(a-u_\lambda)^2}dx.
\end{eqnarray*}

Since the mapping $\lambda\mapsto u_\lambda$ is strictly increasing, we have that
$$\epsilon_\lambda =a(r_n)-u_\lambda(r_n)  \ge a(r_n)-u_{\lambda^*}(r_n):=\varepsilon_n>0$$
and  it infers by $a(x)=\kappa(1-|x|^2)^{\frac23}$,
$$|\Delta a|\le C_n\quad {\rm in}\quad B_{r_n}(0)$$
 for some $C_n>0$, then letting $\theta\in(-2-\sqrt{6},0)$, we have that $2-\frac{\theta^2}{1-2\theta}>0$ and by  H\"{o}lder inequality, we obtain that
\begin{eqnarray*}
&&\lambda(2-\frac{\theta^2}{1-2\theta})\int_{B_{r_n}(0)}\frac1{(a-u_\lambda)^{3-2\theta}}dx
\\&&\le  \int_{B_{r_n}(0)} \frac{4\lambda\epsilon_\lambda^\theta}{(a-u_\lambda)^{3-\theta}}dx+  \frac{\theta^2}{1-2\theta}\int_{B_{r_n}(0)} \frac{C_n}{(a-u_\lambda)^{1-2\theta}}  dx+\frac{\theta^2\varepsilon_n^{2\theta-1}}{1-2\theta} C_n |B_{r_n}(0)|\nonumber\\
\\&&\le   4\lambda^*\varepsilon_n^{\theta} |B_1(0)|^{ \frac{-\theta}{3-2\theta}}\left(\int_{B_{r_n}(0)}\frac{1}{(a-u_\lambda)^{3-2\theta}}dx\right)^{\frac{3-\theta}{3-2\theta}}
\\&&\quad+
C_n |B_1(0)|^{ \frac{2}{3-2\theta}} \left(\int_{B_{r_n}(0)}\frac{1}{(a-u_\lambda)^{3-2\theta}}dx\right)^{\frac{1-2\theta}{3-2\theta}} +\frac{\theta^2\varepsilon_n^{2\theta-1}}{1-2\theta} C_n |B_1(0)|
\end{eqnarray*}
thus,
there exists $c_{44}>0$ independent on $ \lambda$ such that
\begin{equation}\label{4.3.4}
 \int_{B_{r_n}(0)}\frac1{(a-u_\lambda)^{3-2\theta}}dx\le c_{44}.
\end{equation}
When $1\le N\le 7$, $\frac{3N}{2}\le 3-2\theta$ for some $\theta\in(-2-\sqrt{6},0)$, then by Lemma \ref{lm 4.2},
we have that $u_{\lambda}$ has uniformly in $C^{2,\beta}_{loc}(B_1(0))$ with $\beta<\gamma$, then $u_{\lambda^*}$
is a classical solution of (\ref{eq 1.1}) with $\lambda=\lambda^*$ and
$a-u_{\lambda^*}>0$ in $B_1(0)$. \qquad $\Box$

\medskip

\noindent{\bf Proof of Theorem \ref{teo 3}.}
Proposition \ref{pr 4.1} shows that $u_{\lambda^*}$ is a weak solution of (\ref{eq 1.1}) with $\lambda=\lambda^*$.
The stability of $u_\lambda$   follows by Proposition \ref{pr 4.2} for $\lambda\in(0,\lambda_*)$.
 When $\gamma=\frac23$, the stability and regularity of $u_{\lambda^*}$ see Proposition \ref{pr 4.3} and Proposition \ref{pr 4.4}.
 \qquad$\Box$

\bigskip

\noindent{\bf Acknowledgements:}
H. Chen is supported by NSFC, No: 11401270. F. Zhou   is supported by  NSFC, No: 11271133.

\end{document}